\documentclass[12pt]{article}
\usepackage[latin1]{inputenc}
\usepackage{amsmath}
\usepackage{amsfonts}
\usepackage{amssymb}
\usepackage{amsthm}

\usepackage{titlesec}
\titleformat{\subsection}[hang]{\normalfont\bfseries}{\thesubsection}{1em}{}

\usepackage{titlesec}

\titlespacing\section{0pt}{3.5ex plus 0.5ex minus .2ex}{0.3ex plus .2ex}
\titlespacing\subsection{0pt}{2.5ex plus 0.5ex minus .2ex}{0.3ex plus .2ex}
\titlespacing\subsubsection{0pt}{2.5ex plus 0.5ex minus .2ex}{0.3ex plus .2ex}

\usepackage{mathrsfs} 

\usepackage[alphabetic,lite]{amsrefs} 

\usepackage{fullpage}
\usepackage{setspace}
\usepackage{hyperref}
\usepackage{color}
\usepackage{enumerate} 
\usepackage{ulem} 
\usepackage{comment}

\usepackage{marginnote}

\usepackage[all,cmtip]{xy} 

 \usepackage{fancyhdr}
\newtheorem{Thm}{Theorem}

\newtheorem{Lemma}[Thm]{Lemma}

\newtheorem{Cor}[Thm]{Corollary}

\theoremstyle{definition}

\newcommand{\bF}{\mathbb{F}}

\newcommand{\bN}{\mathbb{N}}

\newcommand{\bR}{\mathbb{R}}

\newcommand{\bZ}{\mathbb{Z}}

\newcommand{\cO}{\mathcal{O}}

\newcommand{\ff}{\mathfrak{f}}

\newcommand{\sA}{\mathscr{A}}
\newcommand{\sB}{\mathscr{B}}

\newcommand{\sS}{\mathscr{S}}

\DeclareMathAlphabet{\mathpzc}{OT1}{pzc}{m}{it}

\newcommand{\Proof}{\textbf{Proof.\\}}

\newcommand{\ra}{\rightarrow}

\newcommand{\wt}{\widetilde}

\newcommand{\eps}{\epsilon}

\DeclareMathOperator{\Ind}{Ind}			
\DeclareMathOperator{\cind}{c-ind}  	  
\DeclareMathOperator{\pr}{pr}  	  

\newcommand{\field}{R}

\AtEndDocument{\bigskip{\footnotesize%
		\textsc{University of Cambridge, Cambridge, UK and Duke University, Durham, NC, USA} \par  
		\textit{Mailing address}: Trinity College, Cambridge, CB2 1TQ, UK \par
		\textit{E-mail address}: \texttt{fintzen@maths.cam.ac.uk} and \texttt{fintzen@math.duke.edu}
}}

\begin{document}
\author{Jessica Fintzen}
\title{Supercuspidal representations in non-defining characteristics}
\date{}
\maketitle
\begin{abstract}
 We show that a mod-$\ell$-representation of a $p$-adic group arising from the analogue of Yu's construction is supercuspidal if and only if it arises from a supercuspidal representation of a finite reductive group. This has been previously shown by Henniart and Vigneras under the assumption that the second adjointness holds, a statement that is not yet available in the literature.
\end{abstract}

{
	\renewcommand{\thefootnote}{}  
	\footnotetext{MSC2020: 22E50, 
		20C20 
	} 
	\footnotetext{Keywords: Supercuspidal $\ell$-modular representations, representations of reductive groups over non-archimedean local fields, $p$-adic groups, modular representations}
	\footnotetext{The author was partially supported by NSF Grants DMS-1802234 / DMS-2055230 and DMS-2044643, a Royal Society University Research Fellowship and a Sloan Research Fellowship.}
}



\section{Introduction}
The exhaustive explicit construction and parametrization of supercuspidal irreducible representations of $p$-adic groups with complex coefficients plays a key role in the complex representation theory of $p$-adic groups and beyond. For number theoretic applications it is often desirable to obtain analogous results for representations whose coefficients are valued in an algebraically closed field $\field$ of characteristic $\ell$ different from $p$. In that setting one needs to distinguish between cuspidal and supercuspidal representations. An exhaustive construction of the former, more general notion, is known if the $p$-adic group is tame and $p$ does not divide the order of the Weyl group (\cite{Fi-mod-ell}). This paper concludes the exhaustive construction of supercuspidal irreducible representations for $p$-adic groups in the same setting by determining which of the cuspidal  representations are supercuspidal.

More precisely, we show that if an $\field$-representation arising from the analogue of Yu's construction is supercuspidal, then the representation of a finite reductive group that forms part of the input for the construction has to be supercuspidal as well (Theorem \ref{Thm-main}). Combined with the reverse implication proved by Henniart and Vigneras \cite[Theorem 6.10 and \S6.4.2]{Henniart-Vigneras} and the result that all cuspidal $\field$-representations arise from the analogue of Yu's construction (\cite[Theorem~4.1]{Fi-mod-ell}), we obtain an exhaustive explicit construction of all supercuspidal irreducible $\field$-representations.

Theorem \ref{Thm-main} has previously been proven by Henniart and Vigneras (Theorem~6.10~and~\S6.4.2 in \cite{Henniart-Vigneras}) using different techniques, but only under the assumption that the second adjointness holds in this setting. The second adjointness is so far only proven in the literature for depth-zero representations or if $G$ is a general linear group, a classical group (with $p \neq 2$) or a group of relative rank 1 (\cite{Dat-finitude}). Our approach does not rely on the second adjointness to hold. 

\textbf{Acknowledgments.} The author thanks Guy Henniart and Marie-France Vigneras for helpful discussions related to the topic of this paper.

\section{The main theorem and corollaries}

Let $F$ be a non-archimedean local field of residual characteristic $p$ with ring of integers $\cO$ and residue field $\bF_q$. Let $G$ be a (connected) reductive group over $F$. Let $\field$ be an algebraically closed field of characteristic $\ell$ different from $p$. 
 All representations in this paper are representations with coefficient field $\field$.

For a point $x$ in the enlarged Bruhat--Tits building $\sB(G,F)$ of $G$ over $F$, we write $[x]$ for the image of the point $x$ in the reduced Bruhat--Tits building, 
$G_{[x]}$ for the stabilizer of $[x]$,  $G_{x,0}$ for the parahoric subgroup, and $G_{x,r}$ for the Moy--Prasad filtration subgroup of depth $r$ for a positive real number $r$.

Let $((G_i)_{1 \leq i \leq n+1}, x \in \sB(G_{n+1}, F)\subset \sB(G, F), (r_i)_{1 \leq i \leq n}, \rho, (\phi_i)_{1 \leq i \leq n})$ be an input for the construction of a cuspidal $\field$-representation as in \cite[Section~2.2]{Fi-mod-ell}, i.e. following Yu's construction \cite{Yu} adapted to the mod-$\ell$ setting. If $n>0$, then we assume that $G$ and $G_{n+1}$ split over a tamely ramified extension of $F$ and $p\neq 2$, as in Yu's construction. In the case of $n=0$, the construction recovers the depth-zero representations and we allow $G$ to be wildly ramified and/or $p=2$.
The irreducible $\field$-representation $\rho$ of $(G_{n+1})_{[x]}$ is trivial on $(G_{n+1})_{x,0+}$ and its restriction to $(G_{n+1})_{x,0}$ is a cuspidal representation of $(G_{n+1})_{x,0}/(G_{n+1})_{x,0+}$. Note that the restriction of $\rho$ to $(G_{n+1})_{x,0}$ is semisimple because $(G_{n+1})_{[x]}$ normalizes $(G_{n+1})_{x,0}$. Hence it makes sense to talk about the restriction being supercuspidal or not.

From the tuple $((G_i)_{1 \leq i \leq n+1}, x \in \sB(G_{n+1}, F)\subset \sB(G, F), (r_i)_{1 \leq i \leq n}, \rho, (\phi_i)_{1 \leq i \leq n})$ we obtain a representation $(\wt \rho =\rho \otimes \kappa, V_\rho \otimes V_\kappa)$ of 
$$\wt K=(G_1)_{x,\frac{r_1}{2}}(G_2)_{x,\frac{r_2}{2}}\hdots (G_n)_{x,\frac{r_n}{2}}(G_{n+1})_{[x]}$$
(\cite[Section~2.3]{Fi-mod-ell}), where $\rho$ also denotes the extension of the depth-zero representation $\rho$ from $(G_{n+1})_{[x]}$ to $\wt K$ that is trivial on $(G_1)_{x,\frac{r_1}{2}}(G_2)_{x,\frac{r_2}{2}}\hdots (G_n)_{x,\frac{r_n}{2}}$, such that 
$ \cind_{\wt K}^G(\rho \otimes \kappa)$ is irreducible and cuspidal.

\begin{Thm} \label{Thm-main}
	With the above notation, if $ \cind_{\wt K}^{G(F)}(\rho \otimes \kappa)$ is supercuspidal, then the restriction of $\rho$ to $(G_{n+1})_{x,0}$ is supercuspidal as a representation of $(G_{n+1})_{x,0}/(G_{n+1})_{x,0+}$.
\end{Thm}

This theorem has been proven by Henniart and Vigneras (\cite[Theorem~6.10~and~\S6.4.2]{Henniart-Vigneras}) under the assumption that the second adjointness holds in this setting, which is so far only proven in the literature for depth-zero representations or if $G$ is a general linear group, a classical group (with $p \neq 2$) or a reductive group of relative rank 1 (\cite{Dat-finitude}). 
 
Combining Theorem \ref{Thm-main} with the unconditional result of Henniart and Vigneras  that $ \cind_{\wt K}^G(\rho \otimes \kappa)$ is supercuspidal when the restriction of $\rho$ to $(G_{n+1})_{x,0}$ is supercuspidal as a representation of $(G_{n+1})_{x,0}/(G_{n+1})_{x,0+}$ (\cite[Theorem~6.10~and~\S6.4.2]{Henniart-Vigneras}), we obtain the following unconditional Corollary \ref{Cor-main}.

\begin{Cor} \label{Cor-main}
	With the above notation, $ \cind_{\wt K}^{G(F)}(\rho \otimes \kappa)$ is supercuspidal if and only if the restriction of $\rho$ to $(G_{n+1})_{x,0}$ is supercuspidal as a representation of $(G_{n+1})_{x,0}/(G_{n+1})_{x,0+}$.
\end{Cor}

\Proof Combine Theorem \ref{Thm-main} with Theorem 6.10 and \S6.4.2 of \cite{Henniart-Vigneras}. \qed

Combined with \cite[Theorem~4.1]{Fi-mod-ell} we obtain the following result.

\begin{Cor} \label{Cor-sc}
	Suppose that $G$ splits over a tamely ramified field extension of $F$ and $p$ does not divide the order of the (absolute) Weyl group of $G$. Then all supercuspidal representations of $G(F)$ are of the form $ \cind_{\wt K}^{G(F)}(\rho \otimes \kappa)$ as above where the restriction of $\rho$ to $(G_{n+1})_{x,0}$ is supercuspidal as a representation of $(G_{n+1})_{x,0}/(G_{n+1})_{x,0+}$.
\end{Cor}

\Proof This follows from Corollary \ref{Cor-main} and \cite[Theorem~4.1]{Fi-mod-ell}. \qed

\section{Proof of Theorem \ref{Thm-main}} 
Let $\pi:=\cind_{\wt K}^{G(F)}(\rho \otimes \kappa)$ be a cuspidal irreducible representation as in the previous section.
In this section we will prove that $\pi$ being supercuspidal implies that $\rho$ is supercuspidal, i.e. Theorem \ref{Thm-main}.  This statement is trivially true if $G_{n+1}$ is anisotropic as in this case $(G_{n+1})_{x,0}/(G_{n+1})_{x,0+}$ is also anisotropic and hence all its semisimple representations are supercuspidal. Thus we assume throughout that $G_{n+1}$ is not anisotropic. 

We will eventually prove the desired result by assuming that $\rho$ is not supercuspidal and proving that then $\pi$ is a subquotient of a parabolicly induced representation, i.e. is not supercuspidal either. However, we will not make this assumption until the end of this section to first prove a series of results that hold without this assumption and might be useful on its own for other applications.

Let $T$ be a maximally split, maximal torus of $G_{n+1}$ such that $x$ is contained in the apartment $\sA(T, F)$  of $T$. Let $\lambda$ be a cocharacter of $T$. We write $P_{G_i}(\lambda)$ for the parabolic subgroup of $G_i$ ($1 \leq i \leq n+1$) attached to $\lambda$ as in Section 2.1 and 2.2, in particular Proposition 2.2.9, of \cite{Pseudoreductive2}.  This means $P_{G_i}(\lambda)(F)$ consists of the elements $g \in G_i(F)$ for which the limit of $\lambda(t)g\lambda(t)^{-1}$ as $t$ goes to zero exists (i.e. extends to a map from the affine line to $G_i$). Then the centralizer $Z_{G_i}(\lambda)$ of $\lambda$ is a Levi subgroup of  $P_{G_i}(\lambda)$, which we also denote by  $M_{G_i}(\lambda)$. Let $U_{G_i}(\lambda)$ be the unipotent radical of  $P_{G_i}(\lambda)$ and $\bar U_{G_i}(\lambda)$ the unipotent radical of the opposite parabolic $\bar  P_{G_i}(\lambda)$ of  $P_{G_i}(\lambda)$ with respect to  $M_{G_i}(\lambda)$. 

 Let $\eps>0$ be sufficiently small so that $G_{x,s+} \subset G_{y,s} \subset G_{x,s}$ for $s \in \{ \frac{r_1}{2}, \frac{r_2}{2}, \hdots, \frac{r_n}{2}, 0 \}$ and $y=x +\eps\lambda \in \sA(T, F)$. 
While there is in general no canonical embedding of the Bruhat--Tits building of  $M_{n+1}:=M_{G_{n+1}}(\lambda)$ into the Bruhat--Tits building of $G_{n+1}$, the embedding is unique up to translation by $X_*(Z(M_{n+1})) \otimes_{\bZ} \bR$, where $X_*(Z(M_{n+1}))$ denotes the cocharacters of the center $Z(M_{n+1})$ of $M_{n+1}$, and we fix an embedding of Bruhat--Tits buildings throughout the paper to view $\sB(M_{n+1}, F)$ as a subset of $\sB(G_{n+1}, F)$. We will do the same for all twisted Levi subgroups of $G$ to view all Bruhat--Tits buildings over $F$ as subsets of $\sB(G,F)$.
Then we have  $y=x+\eps\lambda \in \sB(M_{n+1}, F) \subset \sB(G_{n+1}, F) \subset \sB(G,F)$, and
\begin{equation} \label{eq-depthzero}
(G_{n+1})_{y,0}/(G_{n+1})_{y,0+} \simeq (M_{n+1})_{y,0}/(M_{n+1})_{y,0+}  .
\end{equation} 
 
  For $z\in \{x, y\}$ and $s\in \{0, 0+\}$, we set 
$$K_{z,s}=(G_1)_{z,\frac{r_1}{2}}(G_2)_{z,\frac{r_2}{2}}\hdots (G_n)_{z,\frac{r_n}{2}}(G_{n+1})_{z,s}, $$
$$K_{+}=(G_1)_{x,\frac{r_1}{2}+}(G_2)_{x,\frac{r_2}{2}+}\hdots (G_n)_{x,\frac{r_n}{2}+}(G_{n+1})_{x,0+} $$
We might abbreviate the groups $K_{x,0}$ and $K_{x,0+}$ by $K_0$ and $K_{0+}$, respectively, and write $P=P_{G}(\lambda)$, $M=M_{G}(\lambda)$, $U=U_{G}(\lambda)$, $\bar P=\bar P_{G}(\lambda)$ and $\bar U=\bar U_{G}(\lambda)$.
Then $K_{y,0} \subset K_0$ and 
\begin{eqnarray}
K_{y,0} \cap U(F) & = & K_{y,0+} \cap U(F) = K_0 \cap U(F) , \label{eq-1} \\
 K_{y,0} \cap \bar U(F) & = & K_{y,0+} \cap \bar U(F)=  K_{+} \cap \bar U(F)  , \label{eq-2} 
 \\
K_{y,0} &= &(K_{y,0}\cap\bar U(F))(K_{y,0} \cap M(F))(K_{y,0}\cap U(F)) .  \label{eq-Ky0}  \end{eqnarray}

\begin{Lemma}\label{Lemma-decomposed}
	\begin{enumerate}[(a)] 
		\item 	The space of $(K_0 \cap U(F))$-fixed vectors $V_\kappa^{K_0 \cap U(F)}$ of the representation $(\kappa, V_\kappa)$ is non-trivial. 
		\item \label{Lemma-decomposed-b} The representation $(\kappa, V_\kappa)$ is trivial when restricted to the subgroup $K_{y,0}\cap\bar U(F)$.
		\item The subspace $V_\kappa^{K_0 \cap U(F)}$  is preserved under the action of $K_{y,0}$ via $\kappa$.
	\end{enumerate}
\end{Lemma}	
\Proof 
If $\pi$ has depth-zero, $\kappa$ is the trivial one dimensional representation, and hence all statements are trivially true. Thus we may assume $n>0$ and hence that we are in the setting where $G_{n+1}$ splits over a tamely ramified field extension.
Let $E$ be the splitting field of $T$. Since $\lambda$ factors through a maximal split torus, which is contained in a maximal torus the splits over a tamely ramified extension, we may assume without loss of generality that $E$ is tamely ramified over $F$. 
 For $1 \leq i \leq n$, we define 
\begin{eqnarray*} 
	U_i  &=& G(F) \cap \left< U_{\alpha}(E)_{x,\frac{r_i}{2}} \, | \, \alpha \in \Phi(G_i, T) \setminus \Phi(G_{i+1},T), \lambda(\alpha)>0 \right>, \\
		U_{n+1} &=& G(F) \cap \left< U_{\alpha}(E)_{x,0} \, | \, \alpha \in \Phi(G_{n+1}, T), \lambda(\alpha)>0 \right>,
\end{eqnarray*}
where $\Phi(G_i, T)$ denotes the root system of $G_i$ with respect to $T$ over the field $E$ (for $1 \leq i \leq n+1$) and $U_{\alpha}(E)_{x,\frac{r_i}{2}}$ denotes the depth-$\frac{r_i}{2}$ filtration subgroup of the root group $U_{\alpha}(E)$ of $G(E)$ corresponding to $\alpha$ and normalized with respect to the valuation on $E$ that extends the valuation on $F$ used to define the Moy--Prasad filtration.
	Then 
	$$K_0 \cap U(F)=U_1U_2\hdots U_{n+1} .$$
	  Following \cite[Section~2.5]{Fi-Yu-works} we write $V_\kappa = \otimes_{i=1}^n V_{\omega_i}$ so that the action of $\kappa$ restricted to $U_j$ $(1 \leq j \leq n)$ is given by $U_j$ acting on $V_{\omega_k}$ for $k \neq j$ via the character $\hat \phi_k$ defined in \textit{loc. cit.} and on $V_{\omega_j}$ via a Heisenberg representation. The action of $\kappa$ restricted to $U_{n+1}$ arises from $U_{n+1}$ acting on $V_{\omega_k}$ via $\phi_k$ tensored with a composition with a Weil representation, see \textit{loc. cit.} for a precise definition. For $1 \leq j < k \leq n$, the restriction of $\hat \phi_k$ to $U_j$ is trivial by the construction of $\hat \phi_k$. For  $1 \leq k < j \leq n+1$, the restriction of $\hat \phi_k$ to $U_j$ equals the restriction of the character $\phi_k$ from $G_{k+1}(F)$ to $U_j$. Since $U_j$ is contained in the unipotent radical of a parabolic subgroup of $G_{k+1}(F)$, we conclude that the restriction of $\phi_k$ to $U_j$ is trivial (\cites{Tits64, Tits-Kneser-Tits}). Thus 
	$$V_\kappa^{U_1U_2 \hdots U_n} = \bigotimes_{i=1}^n (V_{\omega_i})^{U_i}.$$
	Using the same arguments as in the proof of \cite[Theorem~3.1]{Fi-Yu-works}, we obtain that the space $(V_{\omega_i})^{U_i}$ is nontrivial and that $U_{n+1}$ acts on $(V_{\omega_i})^{U_i}$ via the restriction of the character $\phi_i$ to $U_{n+1}$ for $1 \leq i \leq n$, which we observed above is trivial. Hence 
		$$V_\kappa^{K_0 \cap U(F)} = V_\kappa^{U_1U_2 \hdots U_nU_{n+1}} = \bigotimes_{i=1}^n (V_{\omega_i})^{U_i} \not\simeq \{0 \} .$$
For (b) recall that $\kappa$ restricted to $K_+$ acts via the character $\prod_{1 \leq i \leq n}\hat\phi_i$ (times identity).
 For $1 \leq i \leq n$, we define 
\begin{eqnarray*} 
	\bar U_i^+  &=& G(F) \cap \left< U_{\alpha}(E)_{x,\frac{r_i}{2}+} \, | \, \alpha \in \Phi(G_i, T) \setminus \Phi(G_{i+1},T), \lambda(\alpha)<0 \right>, \\
	\bar U_{n+1}^+ &=& G(F) \cap \left< U_{\alpha}(E)_{x,0+} \, | \, \alpha \in \Phi(G_{n+1}, T), \lambda(\alpha)<0 \right>.
\end{eqnarray*}
Then
$$K_{y,0}\cap\bar U(F)= K_{+} \cap \bar U(F) =\bar U_1^+\bar U_2^+\hdots \bar U_{n+1}^+ .$$
For $1 \leq j \leq i \leq n$, the restriction of $\hat \phi_i$ to $\bar U_j^+$ is trivial by the construction of $\hat \phi_i$ and the definition of $\bar U_j^+$.  For  $1 \leq i < j \leq n+1$, the restriction of $\hat \phi_i$ to $\bar U_j^+$ equals the restriction of the character $\phi_i$ from $G_{i+1}(F)$ to $\bar U_j^+$. Since $\bar U_j^+$ is contained in the unipotent radical of a parabolic subgroup of $G_{k+1}(F)$, the restriction of $\phi_i$ to $\bar U_j^+$ is trivial (\cites{Tits64, Tits-Kneser-Tits}). Hence the restriction of $(\kappa, V_\kappa)$ to $K_{y,0}\cap\bar U(F)= \bar U_1^+\bar U_2^+\hdots \bar U_{n+1}^+$ is trivial.

Claim (c) follows now from Equation  \eqref{eq-Ky0} and the observation that $K_{y,0} \cap M(F)$ normalizes $K_{y,0} \cap U(F)=K_0 \cap U(F)$.

\qed

\begin{Lemma} \label{Lemma-quotient}
	Let $(\rho', V_{\rho'})$ be a representation of $K_{y,0}K_{0+}$ that is trivial on $K_{0+}$. Then there exists a surjection of $\wt K$-representations
	$$\pr:  \cind_{K_{y,0}}^{\wt K}\left(V_{\rho'} \otimes V_\kappa^{K_0\cap U(F)}\right) \twoheadrightarrow \left(\cind_{K_{y,0}K_{0+}}^{\wt K}V_{\rho'}\right)\otimes V_\kappa .$$
\end{Lemma}
\Proof
To ease notation, we abbreviate $ \left(\cind_{K_{y,0}}^{\wt K}\left(\rho' \otimes \kappa \right), \cind_{K_{y,0}}^{\wt K}\left(V_{\rho'} \otimes V_\kappa^{K_0\cap U(F)}\right)\right)$ by $ \left(\sigma_1, V_{\sigma_1}\right)$ and denote 
$\left(\cind_{K_{y,0}K_{0+}}^{\wt K} \rho', \cind_{K_{y,0}K_{0+}}^{\wt K}V_{\rho'}\right)$ by $\left(\sigma_2, V_{\sigma_2}\right)$. Using this notation we need to construct a surjection $\pr$ from $ \left(\sigma_1, V_{\sigma_1}\right)$ to $\left(\sigma_2 \otimes \kappa, V_{\sigma_2} \otimes V_\kappa \right)$.

For $v \in V_{\rho'}, w \in V_\kappa^{K_0\cap U(F)}$, we write $f_{v \otimes w}$ for the element of $V_{\sigma_1}=\cind_{K_{y,0}}^{\wt K}\left(V_{\rho'} \otimes V_\kappa^{K_0\cap U(F)}\right) $ that is supported on $K_{y,0}$ and satisfies $f_{v \otimes w}(1)=v \otimes w$. Then an arbitrary element of $V_{\sigma_1}$ can be written as
$$ \sum_{1 \leq i \leq j} c_i\sigma_1(g_i)f_{v_i \otimes w_i} $$
with $j \in \bN$ and $c_i \in R, g_i \in \wt K, v_i \in V_{\rho'}, w_i \in V_\kappa^{K_0 \cap U(F)}$ for $1 \leq i \leq j$.
For $v \in V_{\rho'}$, we write $f_{v}$ for the element of $V_{\sigma_2}=\cind_{K_{y,0}K_{0+}}^{\wt K}V_{\rho'}$ that is supported on $K_{y,0}K_{0+}$ and satisfies $f_v(1)=v$.
Then we define the morphism $\pr:V_{\sigma_1} \ra V_{\sigma_2} \otimes V_\kappa$ by
$$ \pr\left(\sum_{1 \leq i \leq j} c_i\sigma_1(g_i)f_{v_i \otimes w_i}  \right) = \sum_{1 \leq i \leq j} c_i\left(\sigma_2(g_i)f_{v_i} \otimes \kappa(g_i)w_i \right)  .$$
In order to see that this linear morphism is well defined, we assume that 
$$\sum_{1 \leq i \leq j} c_i\sigma_1(g_i)f_{v_i \otimes w_i} =0$$
and need to show that $\pr(\sum_{1 \leq i \leq j} c_i\sigma_1(g_i)f_{v_i \otimes w_i})=\sum_{1 \leq i \leq j} c_i\left(\sigma_2(g_i)f_{v_i} \otimes \kappa(g_i)w_i  \right)=0$. By considering the support we reduce to the case that $g_i \in K_{y,0}$ for $1 \leq i \leq j$. Let $u_1, u_2, \hdots, u_{j'}$ be a basis for the $R$-linear span of $\{\kappa(g_i)w_i\}_{1 \leq i \leq j}$ and write $\kappa(g_i)w_i=\sum_{1 \leq i' \leq j'}d_{i,i'}u_{i'}$ for some $d_{i,i'}\in R$ for $1 \leq i \leq j$.
Then
\begin{eqnarray*}
0&=&\sum_{1 \leq i \leq j} c_i\sigma_1(g_i)f_{v_i \otimes w_i}(1)=\sum_{1 \leq i \leq j} c_i(\rho' \otimes \kappa)(g_i)(v_i \otimes w_i)=\sum_{1 \leq i \leq j} c_i(\rho'(g_i)v_i \otimes \kappa(g_i)w_i) \\
&=&\sum_{1 \leq i \leq j}\sum_{1 \leq i' \leq j'} c_id_{i,i'}(\rho'(g_i)v_i \otimes u_{i'}) ,
\end{eqnarray*}
which implies 
$$	\sum_{1 \leq i \leq j} c_id_{i,i'}\rho'(g_i)v_i =0 \, \, \text{ for } \, \, 1 \leq i' \leq j' .$$
 Hence
\begin{eqnarray*}
	\sum_{1 \leq i \leq j} c_i\left(\sigma_2(g_i)f_{v_i} \otimes \kappa(g_i)w_i \right)
	&=& \sum_{1 \leq i \leq j}\sum_{1 \leq i' \leq j'}  c_i\left(f_{\rho'(g_i)v_i} \otimes d_{i,i'}u_{i'} \right) \\
	&=& \sum_{1 \leq i' \leq j'}\left(\left(\sum_{1 \leq i \leq j}  c_id_{i,i'}f_{\rho'(g_i)v_i}\right) \otimes u_{i'} \right)
=0 .
\end{eqnarray*}
Therefore $\pr$ is well defined and is by construction a $\wt K$-homomorphism. It remains to show that the map $\pr: V_{\sigma_1}\ra V_{\sigma_2} \otimes V_\kappa$ is surjective.
Let $v$ be a non-zero element of $V_{\sigma_2} \otimes V_\kappa= \left(\cind_{K_{y,0}K_{0+}}^{\wt K}V_{\rho'}\right)\otimes V_\kappa$. Then there exists a positive integer $j$, and elements $v_i \in V_{\rho'}, w_i \in V_\kappa$ and $g_i \in \wt K$ for $1 \leq i \leq j$ such that 
$$v=\sum_{1 \leq i \leq j} \sigma_2(g_i)f_{v_i} \otimes w_i .$$
Since the restriction of $\kappa$ to $K_{0+}$ is irreducible, there exist an element $w \in V_\kappa^{K_0 \cap U(F)}$, an integer $j'$, elements $h_{i'} \in K_{0+}$ and	$c_{i,{i'}} \in R$ for $1 \leq i \leq j$ and $1 \leq {i'} \leq j'$ such that 
$$\kappa(g_i^{-1})(w_i)=\sum_{1 \leq {i'} \leq j'} c_{i,{i'}} \kappa(h_{i'})(w) .$$	
Thus, using that $(\rho', V_{\rho'})$ is trivial on $K_{0+}$, we obtain
\begin{eqnarray*}
	v&=& \sum_{1 \leq i \leq j} \sigma_2(g_i)f_{v_i} \otimes w_i \\
	&=& \sum_{1 \leq i \leq j} (\sigma_2 \otimes \kappa)(g_i)(f_{v_i} \otimes \kappa(g_i^{-1})(w_i)) \\
	&=& \sum_{1 \leq i \leq j}\sum_{1 \leq {i'} \leq j'}  (\sigma_2 \otimes \kappa)(g_i)(f_{v_i} \otimes  c_{i,{i'}} \kappa(h_{i'})(w)) \\
	&=& \sum_{1 \leq i \leq j}\sum_{1 \leq {i'} \leq j'}  c_{i,{i'}} (\sigma_2 \otimes \kappa)(g_ih_{i'})(f_{v_i} \otimes  w) \\
	&=& \pr\left(\sum_{1 \leq i \leq j}\sum_{1 \leq {i'} \leq j'}  c_{i,{i'}} \sigma_1(g_ih_{i'})f_{v_i \otimes  w}\right) 	.
\end{eqnarray*}
Hence $\pr$ is surjective. \qed

This allows us to prove the following key lemma for the proof of Theorem \ref{Thm-main}.

\begin{Lemma}\label{Lemma-subquotient}
	If $\rho$ is not supercuspidal, then there exists a maximally split, maximal torus $T$ of $G_{n+1}$ whose apartment contains $x$, a cocharacter $\lambda$ of $T$ and a representation $(\rho', V_{\rho'})$ of $K_{y,0}$ (with $y=x+\eps \lambda$ as above) that is trivial on $K_{y,0+}$ such that the representation $(\rho \otimes \kappa, V_\rho \otimes V_\kappa)$ is a subquotient of $\cind^{\wt K}_{K_{y,0}} (V_{\rho'} \otimes V_\kappa^{K_0 \cap U(F)})$. 
	
	The cocharacter $\lambda$ can be chosen so that $M:=M_{G}(\lambda)$ is the centralizer of the maximal split torus in the center of $M_{n+1}:=M_{G_{n+1}}(\lambda)$ and $\eps>0$ can be chosen so that the point $y=x+\eps\lambda \in \sB(M_{n+1},F)\subset \sB(G,F)$ is contained in a facet of minimal dimension of $\sB(M_{n+1},F)$ and
	\begin{equation}
	\label{eq-generic-embedding}
	 \sum_{i=1}^n\left(\dim\left((G_i)_{y,\frac{r_i}{2}}/(G_i)_{y,r\frac{r_i}{2}+}\right)-\dim\left((M_{G_i}(\lambda))_{y,\frac{r_i}{2}}/(M_{G_i}(\lambda))_{y,\frac{r_i}{2}+}\right)\right)=0 .
	\end{equation}
\end{Lemma}
\Proof
 Suppose $\rho$ is not supercuspidal. Recall that $\rho|_{(G_{n+1})_{x,0}}$ is semisimple as $(G_{n+1})_{x,0}$ is normal inside $(G_{n+1})_{[x]}$. Let $\rho_1$ be an irreducible quotient of $\rho|_{(G_{n+1})_{x,0}}$, viewed as a representation of $(G_{n+1})_{x,0}/(G_{n+1})_{x,0+}$. We denote by $\mathbf G$ the connected reductive group over $\bF_q$ that satisfies for any unramified field extension $E$ of $F$ with residue field $\ff_E$ that $\mathbf G(\ff_E)=G_{n+1}(E)_{x,0}/G_{n+1}(E)_{x,0+}$. Let $\mathbf P$ be a proper parabolic subgroup of $\mathbf G$ with Levi subgroup $\mathbf M$, and $\rho'$ a representation of $\mathbf M(\bF_q)$ such that $\rho_1$ is a subquotient of the parabolic induction $\Ind_{\mathbf P(\bF_q)}^{(G_{n+1})_{x,0}/(G_{n+1})_{x,0+}}\rho'$. Let $\mathbf S$ be a maximal split torus of $\mathbf M$ and $\sS$ the split torus defined over $\cO$ contained in the parahoric group scheme attached to $G_{n+1}$ and $x$ such that the special fiber of $\sS$ is $\mathbf S$. We denote the generic fiber $\sS_F$ of $\sS$ by $S$. Note that $S$ is a maximal split torus of $G_{n+1}$.  Let $\mathcal C$ be the  split subtorus of $\mathcal S$ whose special fiber $\mathbf C := \mathcal C_{\bF_q}$ is the maximal split torus in the center of $\mathbf M$. Let $M$ be the centralizer of $C:=\mathcal C_F$ in $G$. Then $M$ is a Levi subgroup of $G$ and there exists a cocharacter  $\lambda \in X_*(S)$ such that $M=M_G(\lambda)$ (e.g. by \cite[Proposition~2.2.9]{Pseudoreductive2} combined with the fact that Levi subgroups of a fixed parabolic are rationally conjugate).
Choosing a maximally split, maximal torus $T$ of $G_{n+1}$ containing $S$, we can perform the above constructions to obtain a parabolic subgroup $P_{G_{i}}(\lambda)$ of $G_i$ ($1 \leq i \leq n+1$) with Levi subgroup $M_i:=M_{G_i}(\lambda)=Z_{G_i}(\lambda)$ and a point $y=x +\eps\lambda$ in the apartment $\sA(T, F)$. Note that $M_{n+1}$ is the centralizer of $C$ in $G_{n+1}$, because $M$ is the centralizer of $C$ in $G$. Hence by Equation \eqref{eq-depthzero} and \cite[Proposition~6.4(1)]{MP2}, the point $y$ is a minimal facet of the building $\sB(M_{n+1},F)$. Moreover, since $\mathbf C$ is the maximal split torus in the center of $\mathbf M$ and $\mathbf M(\bF_q)=(M_{n+1})_{y,0}/(M_{n+1})_{y,0+}$, the torus $C$ is the maximal split torus in the center of $M_{n+1}$. Hence $M$ is the centralizer of the maximal split torus in the center of $M_{n+1}$, as desired. Moreover, by the definition of $M_{G_i}(\lambda)$, Equation \eqref{eq-generic-embedding} is satisfied by all but finitely many $\eps$ in the open interval $(0,1)$. Hence we may choose $\eps >0$ such that Equation \eqref{eq-generic-embedding} holds true.

Since we have
$$\mathbf M(\bF_q)=K_{y,0}/K_{y,0+}=K_{y,0}K_{0+}/K_{y,0+}K_{0+} ,$$
we may view $\rho'$ as a representation of $K_{y,0}K_{0+}$ via inflation. Note that the image of $K_{y,0}K_{0+}$ in $K_{x,0}/K_{0+}\simeq (G_{n+1})_{x,0}/(G_{n+1})_{x,0+}$ is $\mathbf P(\bF_q)$.
 Viewing $\rho$ and $\rho_1$ as representations of $\wt K$ and $K_{x,0}$, respectively, by asking them to be trivial on $(G_1)_{x,\frac{r_1}{2}}(G_2)_{x,\frac{r_2}{2}}\hdots (G_n)_{x,\frac{r_n}{2}}$, we have by Frobenius reciprocity that the irreducible representation  $\rho$ is a quotient of $\cind_{K_{x,0}}^{\wt K} \rho_1$ and therefore a subquotient of 
 $$\cind_{K_{x,0}}^{\wt K} \cind_{K_{y,0}K_{0+}}^{K_{x,0}} \rho' =\cind_{K_{y,0}K_{0+}}^{\wt K} \rho' .$$ 
 Therefore $(\rho \otimes \kappa, V_\rho \otimes V_\kappa)$ is a subquotient of $((\cind_{K_{y,0}K_{0+}}^{\wt K} \rho') \otimes \kappa, (\cind_{K_{y,0}K_{0+}}^{\wt K} V_\rho') \otimes V_\kappa)$. From Lemma \ref{Lemma-quotient}, we deduce that $(\rho \otimes \kappa, V_\rho \otimes V_\kappa)$ is a subquotient of $\cind^{\wt K}_{K_{y,0}} (V_{\rho'} \otimes V_\kappa^{K_0 \cap U(F)})$.
\qed

\textbf{Proof of Theorem \ref{Thm-main}.}\\
 Suppose $\rho$ is not supercuspidal. We need to prove that  $\pi:=\cind_{\wt K}^{G(F)}(\rho \otimes \kappa)$ is not supercuspidal. We let $\lambda$ be as given by Lemma \ref{Lemma-subquotient}, which provides us with a point $y=x+\eps \lambda$ and a parabolic subgroup $P=P_G(\lambda)$ of $G$ with Levi $M=M_G(\lambda)$ and unipotent radical $U$ as above. Then the representation $(\rho \otimes \kappa, V_\rho \otimes V_\kappa)$ is a subquotient of $\cind^{\wt K}_{K_{y,0}} (V_{\rho'} \otimes V_\kappa^{K_0 \cap U(F)})$.	Hence $\pi$ is a subquotient of
$\cind^{G(F)}_{K_{y,0}} (V_{\rho'} \otimes V_\kappa^{K_0 \cap U(F)})$. We will show that
the latter is isomorphic to a parabolic induction of a smooth representation from $P(F)$, 
which will imply that $\pi$ is not supercuspidal and hence finish the proof.

Recall from Equations \eqref{eq-1}, \eqref{eq-2} and \eqref{eq-Ky0} that 
$$K_{y,0} =(K_{y,0}\cap\bar U(F))(K_{y,0} \cap M(F))(K_{y,0}\cap U(F)) $$
and that $K_{y,0}\cap\bar U(F) =K_{y,0+}\cap\bar U(F) $ and $K_{y,0}\cap U(F)=K_{y,0+}\cap U(F)$.
Moreover, by Lemma \ref{Lemma-decomposed}(\ref{Lemma-decomposed-b}) and since $K_{0}\supset K_{y,0}$ and $(\rho', V_{\rho'})$ is trivial on $K_{y,0+}$, the restriction of $V_{\rho'} \otimes V_\kappa^{K_0 \cap U(F)}$ to $K_{y,0}\cap\bar U(F)$ and to $K_{y,0}\cap U(F)$ is trivial. Hence the pair 
$$(K_{y,0}, (\rho'\otimes \kappa, V_{\rho'} \otimes V_\kappa^{K_0 \cap U(F)}))$$ 
is decomposed over the pair 
$$(K_{y,0}\cap M(F), ((\rho'\otimes \kappa)|_{K_{y,0}\cap M(F)}, V_{\rho'} \otimes V_\kappa^{K_0 \cap U(F)}))$$ 
with respect to $\bar P$ as in the notation of \cite[p.~245]{Blondel04}.

We write  
$$K_{y,+}=(G_1)_{y,\frac{r_1}{2}+}(G_2)_{y,\frac{r_2}{2}+}\hdots (G_n)_{y,\frac{r_n}{2}+}(G_{n+1})_{y,0+} $$
and note that the action of $K_{y,+}$ on $V_{\rho'} \otimes V_\kappa^{K_0 \cap U(F)}$  via $\rho'\otimes \kappa$ is given by $\prod_{1 \leq i \leq n}\hat\phi_i$ (times identity). Let $(\pi',V')$ be an irreducible smooth representation of $G(F)$. Then we write $V'^{(K_{y,+}, \prod \hat \phi_i)}$ for the subspace of $V'$ on which $K_{y,+}$ acts via $\prod_{1 \leq i \leq n}\hat\phi_i$. Since $y$ is contained in a facet of minimal dimension of $\sB(M_{n+1},F)$ and Equation \eqref{eq-generic-embedding} holds by Lemma \ref{Lemma-subquotient} (which ensures that the embedding of the Bruhat--Tits buildings is $(0, \frac{r_n}{2}, \hdots, \frac{r_1}{2})$-generic relative to $y$  as defined by \cite[3.5~Definition]{Kim-Yu}, see also \cite[p.~341]{Fi-exhaustion}) we can apply the proof of \cite[6.3~Theorem]{Kim-Yu} to obtain that the restriction of the Jacquet functor with respect to $\bar U$ to the subspace $V'^{(K_{y,+}, \prod \hat \phi_i)}$ is injective. Note that while Kim and Yu work with complex coefficients in \cite{Kim-Yu}, their proof and \cite[Proposition~6.7]{MP2}, on which the proof relies, also work with coefficients in the field $\field$. Hence the Jacquet functor with respect to $\bar U$ is also injective when restricted  to the maximal subspace of $V'$ that is isomorphic to a direct sum of copies of $(\rho'\otimes \kappa, V_{\rho'} \otimes V_\kappa^{K_0 \cap U(F)})$ as a $K_{y,0}$-representation.
Therefore, the pair
$$(K_{y,0}, (\rho'\otimes \kappa, V_{\rho'} \otimes V_\kappa^{K_0 \cap U(F)}))$$ 
is a cover of 
$$(K_{y,0}\cap M(F), ((\rho'\otimes \kappa)|_{K_{y,0}\cap M(F)}, V_{\rho'} \otimes V_\kappa^{K_0 \cap U(F)}))$$ 
 with respect to $\bar P$ as in \cite[p.~246 and Corollaire de Proposition~2]{Blondel04}. Thus by \cite[Théorème~2]{Blondel04} we have an isomorphism of $G(F)$-representations
 $$ \cind_{K_{y,0}}^{G(F)}(V_{\rho'} \otimes V_\kappa^{K_0 \cap U(F)}) \simeq \cind_{(K_{y,0}\cap M(F))U(F)}^{G(F)}(V_{\rho'} \otimes V_\kappa^{K_0 \cap U(F)}), $$
 where $U(F)$ acts trivially on $ V_{\rho'} \otimes V_\kappa^{K_0 \cap U(F)}$.
 Therefore $\pi$ is a subquotient of
 \begin{eqnarray*}
 \cind_{(K_{y,0}\cap M(F))U(F)}^{G(F)}  (V_{\rho'} \otimes V_\kappa^{K_0 \cap U(F)}) 
& \simeq & \cind_{P(F)}^{G(F)}\cind^{M(F)U(F)}_{(K_{y,0}\cap M(F))U(F)} (V_{\rho'} \otimes V_\kappa^{K_0 \cap U(F)}) \\
 & \simeq & \Ind_{P(F)}^{G(F)}\left(\cind^{M(F)}_{K_{y,0}\cap M(F)}(V_{\rho'} \otimes V_\kappa^{K_0 \cap U(F)})\right)   ,
 \end{eqnarray*}
 where $\Ind_{P(F)}^{G(F)}$ denotes the (unnormalized) parabolic induction. This is a contradiction to $\pi$ being supercuspidal. \qed

\bibliography{Fintzenbib}

\begin{bibdiv}
\begin{biblist}

\bib{Blondel04}{article}{
      author={Blondel, Corinne},
       title={Quelques propri\'{e}t\'{e}s des paires couvrantes},
        date={2005},
        ISSN={0025-5831},
     journal={Math. Ann.},
      volume={331},
      number={2},
       pages={243\ndash 257},
         url={https://doi.org/10.1007/s00208-004-0579-1},
      review={\MR{2115455}},
}

\bib{Pseudoreductive2}{book}{
      author={Conrad, Brian},
      author={Gabber, Ofer},
      author={Prasad, Gopal},
       title={Pseudo-reductive groups},
     edition={Second},
      series={New Mathematical Monographs},
   publisher={Cambridge University Press, Cambridge},
        date={2015},
      volume={26},
        ISBN={978-1-107-08723-1},
         url={https://doi.org/10.1017/CBO9781316092439},
      review={\MR{3362817}},
}

\bib{Dat-finitude}{article}{
      author={Dat, Jean-Francois},
       title={Finitude pour les repr\'{e}sentations lisses de groupes
  {$p$}-adiques},
        date={2009},
        ISSN={1474-7480},
     journal={J. Inst. Math. Jussieu},
      volume={8},
      number={2},
       pages={261\ndash 333},
         url={https://doi.org/10.1017/S1474748008000054},
      review={\MR{2485794}},
}

\bib{Fi-Yu-works}{article}{
      author={Fintzen, Jessica},
       title={On the construction of tame supercuspidal representations},
        date={2021},
        ISSN={0010-437X},
     journal={Compos. Math.},
      volume={157},
      number={12},
       pages={2733\ndash 2746},
         url={https://doi.org/10.1112/S0010437X21007636},
      review={\MR{4357723}},
}

\bib{Fi-exhaustion}{article}{
      author={Fintzen, Jessica},
       title={Types for tame {$p$}-adic groups},
        date={2021},
        ISSN={0003-486X},
     journal={Ann. of Math. (2)},
      volume={193},
      number={1},
       pages={303\ndash 346},
         url={https://doi.org/10.4007/annals.2021.193.1.4},
      review={\MR{4199732}},
}

\bib{Fi-mod-ell}{misc}{
      author={Fintzen, Jessica},
       title={Tame cuspidal representations in non-defining characteristics},
        note={Preprint, available at
  \url{https://arxiv.org/pdf/1905.06374.pdf}.},
}

\bib{Henniart-Vigneras}{misc}{
      author={Henniart, Guy},
      author={Vigneras, Marie-France},
       title={Representations of a reductive p-adic group in characteristic
  distinct from p},
        note={To appear in {T}unisian {J}. of {M}athematics, available at
  \url{https://arxiv.org/pdf/2010.06462v2.pdf}},
}

\bib{Kim-Yu}{incollection}{
      author={Kim, Ju-Lee},
      author={Yu, Jiu-Kang},
       title={Construction of tame types},
        date={2017},
   booktitle={Representation theory, number theory, and invariant theory},
      series={Progr. Math.},
      volume={323},
   publisher={Birkh\"auser/Springer, Cham},
       pages={337\ndash 357},
      review={\MR{3753917}},
}

\bib{MP2}{article}{
      author={Moy, Allen},
      author={Prasad, Gopal},
       title={Jacquet functors and unrefined minimal {$K$}-types},
        date={1996},
        ISSN={0010-2571},
     journal={Comment. Math. Helv.},
      volume={71},
      number={1},
       pages={98\ndash 121},
}

\bib{Tits64}{article}{
      author={Tits, J.},
       title={Algebraic and abstract simple groups},
        date={1964},
        ISSN={0003-486X},
     journal={Ann. of Math. (2)},
      volume={80},
       pages={313\ndash 329},
         url={https://doi.org/10.2307/1970394},
      review={\MR{164968}},
}

\bib{Tits-Kneser-Tits}{incollection}{
      author={Tits, Jacques},
       title={Groupes de {W}hitehead de groupes alg\'{e}briques simples sur un
  corps (d'apr\`es {V}. {P}. {P}latonov et al.)},
        date={1978},
   booktitle={S\'{e}minaire {B}ourbaki, 29e ann\'{e}e (1976/77)},
      series={Lecture Notes in Math.},
      volume={677},
   publisher={Springer, Berlin},
       pages={Exp. No. 505, pp. 218\ndash 236},
      review={\MR{521771}},
}

\bib{Yu}{article}{
      author={Yu, Jiu-Kang},
       title={Construction of tame supercuspidal representations},
        date={2001},
        ISSN={0894-0347},
     journal={J. Amer. Math. Soc.},
      volume={14},
      number={3},
       pages={579\ndash 622 (electronic)},
}

\end{biblist}
\end{bibdiv}

\end{document}